\pdfoutput=1
\documentclass[twoside,12pt,a4paper]{cmft}

\pagespan{1}{?}

\theoremstyle{remark}

\newtheorem*{acknowledgement}{Acknowledgment}

\usepackage{graphicx}
\usepackage{hyperref}
\usepackage{caption}
\usepackage{subcaption}
\usepackage{mathtools}
\newcommand{\defeq}{\vcentcolon=}

\usepackage{tikz}
\usetikzlibrary{shapes,arrows}

\newcounter{minutes}
\divide\time by 60
\newcounter{hours}
\multiply\time by 60
\addtocounter{minutes}{-\time}

\begin{document}

\title[Conformal map and harmonic measure of the Bunimovich stadium]{Conformal map and harmonic measure 
\\ of the Bunimovich stadium}
\date{\today}

\author[V. K. Varma]{Vipin Kerala Varma}
\email{vvarma@ictp.it}
\address{The \textit{Abdus Salam} International Centre for Theoretical Physics,
         Condensed Matter and Statistical Physics Section,
         Strada Costiera 11,
         IT-34151 Trieste,
         Italy}


\begin{abstract}
\noindent We consider the conformal mapping of the Bunimovich stadium, a region enclosed by a Jordan curve with four smooth corners, primarily 
in the context of 
a particle undergoing Brownian motion within its closed geometry with Dirichlet boundary conditions. 
A Chebyshev weighting of the solutions of Symm's 
integral equation is employed to give a numerical conformal map of the region onto the canonical domain of the unit disk 
in the 
complex plane. As a measure of the accuracy of the transformation, the domes' harmonic measure evaluated at the centre of the stadium is 
thereby extracted and is compared with results obtained from
Schwarz-Christoffel transformations and Monte Carlo simulations; the pros and cons of the method are reiterated.
\end{abstract}

\maketitle

\section{Introduction}
Conformal maps are transformations between domains of the complex plane, with local angles being preserved. This 
latter constraint restricts the domains of the complex plane where such a map may be explicitly computed; 
in all other cases maps which are only approximately conformal may be evaluated, whose accuracy may be increased 
either with the development of better algorithms or with higher precision computer implementations of 
existing algorithms.\\
Indeed both exact and approximate maps have long been computed for a variety of scenarios, for instance in 
classic problems of mathematics to determine
conformally invariant solutions of Laplace's equation for heat flow \cite{SL}; in industrial applications for design of 
airfoils and porthole windows \cite{DT}; in civil engineering for determining vibration modes of clamped structures and 
plates \cite{SL}; 
in maritime navigation for map-making \cite{Porter}; in physics for describing field theories of critical phenomena in a 
range of 
two-dimensional statistical models \cite{Belavin}, to name a select few. 

A primary utility of such transformations lie in studies of two dimensional systems confined to particular geometries 
on the plane. This is because, by virtue of the Riemann mapping theorem, a region of the complex plane may 
be conformally 
mapped into 
the interior of the unit disk or any other \textit{canonical} domain; the reformulated problem on the canonical domain 
often can be solved 
more straightforwardly. Moreover, the continuous extension of this mapping to the boundary of the regions is established by the 
Carath\'{e}odory-Osgood theorem \cite{Hen3}; indeed the reverse is also true: the boundary mapping may be extended into 
the interior of the region, 
which reduces the complexity of constructing 
the initial map from a two-dimensional problem to a one-dimensional problem \cite{Porter}. However while the Riemann 
mapping theorem 
guarantees the existence of such a transformation, computing the actual map can be a difficult matter depending on the 
geometry under consideration; the first such 
constructions were by Christoffel and Schwarz \cite{DT}, which may be taken to constitute a ``constructive'' 
proof of the Riemann mapping theorem.

The paper is organised as follows: in Section \ref{sec: Bunimovich} we describe the problem 
of Brownian motion, details of the Bunimovich geometry, and the questions addressed; 
in Section \ref{sec: MC} we simulate the Brownian motion in the stadium stochastically using the Monte Carlo algorithm and extract crude estimates of certain 
probability measures.
In the next two Sections \ref{sec: SC-map} and \ref{sec: Symms} we reformulate 
the problem into the deterministic language of conformal transformations, and perform the conformal mapping first by 
the Schwarz-Christoffel method and then by solving Symm's integral equation via Chebyshev weighted solutions.
We conclude 
in Section \ref{sec: Conclusions} commenting on the accuracies of the various solutions and highlighting the differences 
of Brownian motion between the stadium and the rectangular geometry.

\section{Brownian motion in the Bunimovich stadium} \label{sec: Bunimovich}
Brownian motion is a paradigmatic example of a stochastic process in which a massive particle frequently 
collides with the many surrounding lighter particles, such that a straightforward application of 
Newtonian mechanics breaks down. The dependence of the process's dynamics and properties
on the confining geometries has been investigated in the context of fluid transport \cite{Haenggi} and as a 
problem in multiprecision computing \cite{SIAM}. In the latter work, the probability of a Brownian particle 
in a rectangular geometry to hit the sides with length $1/L$ (the other pair of sides being of unit length), 
having started from the centre, was shown to be exactly given by $p = \cfrac{2}{\pi} \arcsin{k_{L^2}}$, where 
$k$ are singular moduli of elliptic integrals. An analytic treatment of the Brownian process through conformal maps of the geometry 
was also undertaken in the same work, which in turn gave a highly precise estimate for $p$.\\
These interesting results and observations led us to wonder what the situation might turn out to be for more complex geometries. 
A simple extension (with maximal symmetries preserved)
we consider is to attach smooth edges in the form of semi-circles to the two ends of a rectangle, and then pose 
the question, similar as in Ref. \cite{SIAM}: 
with what accuracy can a conformal map of the region be determined through readily programmable techniques of conformal 
mapping, which in turn determines the accuracy of $p$?\\
This Brownian motion problem of determining $p$, when formulated deterministally \cite{SIAM}, is equivalent to determining 
the harmonic measure of $\Gamma_d$ with respect to $\Omega$ at the centre subject to the Dirichlet boundary conditions
\begin{equation}
\label{eq: Poisson}
 \Delta u = 0 \hspace{2em} \textrm{on} \hspace{0.5em}\Omega , \hspace{2em} u|_{\Gamma_d} = 1, \hspace{1em} u|_{\Gamma_s} = 0,
\end{equation}
where $\Delta$ is the discrete Laplacian; the required probability or harmonic measure $p = u(0,0)$. It might certainly be possible to 
obtain better estimates for $p$ using a direct solver of Laplace's equation \eqref{eq: Poisson} that utilise finite element methods, 
which we leave for future work. Our motivation, however, primarily in Section \ref{sec: CMap}, is to use the convergence rate 
of the harmonic measure as a gauge of the accuracy of the stadium's conformal map.
We also note that the stadium has the property of being one among the simplest geometries where 
quantum chaos is exhibited by a Schr\"{o}dinger particle \cite{Gutzwiller}, which is (at least formally) equivalent to a Brownian particle 
with an imaginary diffusion constant.\\
\begin{figure}[h!]
\centering
\includegraphics[width=10.5cm]{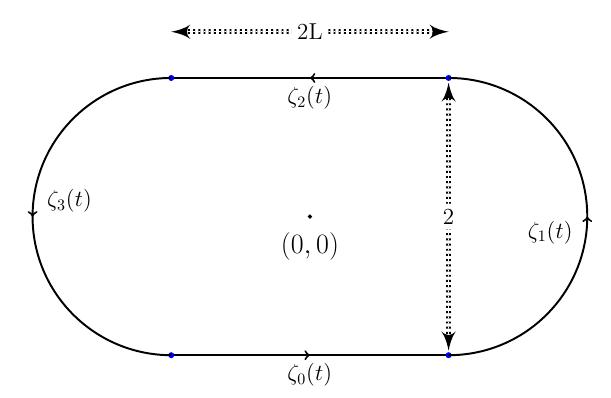}
\caption{Bunimovich stadium $\Omega$ enclosed by the four analytical boundaries $\Gamma \defeq \cup_{i=1}^4 \zeta_i (t)$; 
the parameter $t$ increases from -1 to 1 in each curve in the counter-clockwise direction. The ``smooth'' corners and the 
image of the conformal origin are indicated.}
\label{fig: stadium}
\end{figure}
The Bunimovich stadium, of side lengths $2L$ and domes of unit radii at the ends, that we study is sketched in Fig. \ref{fig: stadium}. 
In this case, as opposed to rectangular geometries, 
one needs to numerically evaluate the conformal map to a given canonical domain; our domain of 
choice will be the interior of the unit disk in the complex $w-$plane
\begin{equation}
 C = \{w: |w| < 1\}. 
\end{equation}
In order to consider the upper half plane $\Re{(w)} > 0$ as the canonical domain, we may first map onto $C$; the conformal 
transformation from the upper half plane to the unit disk is readily effected by the M\"{o}bius transformation 
$f(w) = i\cfrac{w-i}{w+i}$, where $i^2 = 1$.

The four Jordan curves of the Bunimovich stadium of side length $L$ and domes of unit radii are parametrised as 
\begin{eqnarray}
\label{eq: Bunimovich}
 \zeta_0 (t) &\defeq& Lt - i, \hspace{2em}\zeta_1 (t) \defeq L + \sqrt{1 - t^2} + it,\nonumber \\
 \zeta_2 (t) &\defeq& -\zeta_0 (t), \hspace{2em}\zeta_3 (t) \defeq -\zeta_1 (t),
\end{eqnarray}
where $-1 \le t \le 1$ for each arc $\zeta_i (t)$ proceeding in the counter-clockwise direction as indicated by 
the arrows in Fig. \ref{fig: stadium}. The curve $\Gamma_d$ whose harmonic measure $p$ we wish to determine is given by 
$\Gamma_d \defeq \zeta_1 \cup \zeta_3$; clearly the harmonic measure of $\Gamma_s \defeq \zeta_0 \cup \zeta_2$ at 
the centre of the stadium is $1-p$.\\
We will focus on the harmonic measure of $\Gamma_d$ for the case $L=1$ in order to test the accuracy of the map i.e. when the unit-radii 
domes are attached to a square with its corresponding ends removed. A range of $L$ values will also be investigated to analyse the 
effect of stretching or contracting the stadium in the horizontal direction, which will manifest additional contrasts in the 
behaviour of Brownian motion between the rectangular and stadium geometries.

\section{Monte Carlo sampling} \label{sec: MC}
The Monte Carlo algorithm is well suited to simulate processes that are random or stochastic in nature. 
This serves to give a first estimate of $\Gamma_d$'s harmonic measure.
The Monte Carlo procedure we follow is identical to that in Ref. \cite{SIAM} employed for the rectangular geometry, and 
is outlined below: 
if the motion is isotropic, then at any given point $P$ in space the particle 
has equal probability, after a certain time, of reaching any point on a circle drawn with $P$ as the centre. \\
Therefore at 
any point $P$, starting from the centre, we draw the largest circle that fits into $\Omega$ and randomly place the particle on any 
point on this circle. The process is repeated until the particle's distance $r_i$ from the arc $\zeta_i$ is less than or equal to a 
certain specified small value $h \ll 1$; upon which, we assert that the particle has hit the arc $\zeta_i$. This is repeated over 
$N$ trials, initialising a random seed at the commencement for the Monte Carlo runs, 
to give a statistical estimate of the harmonic measure of $\Gamma_d$ at the centre by 
counting all the hits. A counter is updated if the hit is on $\Gamma_d$ else we move onto the next trial; the required probability $p$ is 
approximated by 
\begin{equation}
 \label{eq: hits}
 p \approx \frac{\textrm{(\# hits)}_{\Gamma_d}}{N}.
\end{equation}

\begin{center}
    \begin{table}[h!]
    \begin{tabular}{| l | l | l |}
    \hline
    $N$ & $h$ & $p \approx \textrm{hits}/N$ \\ \hline
    $10^7$ & $10^{-3}$ & 0.281814  \\ 
     & $10^{-5}$ & 0.281852 \\ \hline
    $10^8$ & $10^{-3}$ & 0.281647  \\ 
     & $10^{-5}$ & 0.281727 \\ \hline
    $10^9$ & $10^{-3}$ & 0.281686 \\ 
     & $10^{-5}$ &  0.281802\\ 
    \hline
    \end{tabular}
     \caption{Statistical estimates of the harmonic measure of $\Gamma_d$ at the centre of the stadium from $N$ Monte Carlo runs.}
     \label{tab: MC}
    \end{table}
\end{center}
In Table \ref{tab: MC} we have listed the results of our Monte Carlo simulations for $N = 10^7, 10^8, 10^9$ runs and $h = 10^{-5}, 10^{-3}$. 
Clearly the statistical nature of the errors introduced allow us to be confident of only the first 3 digits i.e. $p \approx 0.281$; this level 
of uncertainty, for the number of runs $N$ considered, is also substantiated by a staightforward theory of statistics \cite{SIAM} which 
estimates that the absolute error in each of the calculated $p$ values can be shown to be around $10^{-4}$.

\section{Conformal mapping} \label{sec: CMap}
In the previous section we treated the Brownian problem in stochastic terms simulated by Monte Carlo runs; 
the accuracy in determining the harmonic measure was limited by the statistical nature of the procedure. 
In this section we construct conformal maps of the stadium $\Omega$ (and regions closely approximating it), and the boundary mapping $\Gamma$ on the 
circumference of $C$. This establishes the locations of the prevertices (i.e. images of the stadium's vertices), from which the accuracy of the map and the harmonic measure $p$ may 
be inferred. The latter is possible because, by observing 
that conformal transplants of harmonic functions such as $u$ preserve the harmonic measure \cite{SIAM}, we may map the original geometry onto the 
canonical domain $C$ and solve the problem much more easily in $C$; the difficulty is, of course, in constructing an 
accurate conformal transformation. \\
For the actual computation of the Brownian particle's probability to hit the domes, we utilise the important fact \cite[Lemma 10.2]{SIAM} 
that if the arc (or union of arcs) $\tilde{\Gamma}$ of the unit circle has length $2\pi p$, the harmonic measure of $\tilde{\Gamma}$ with 
respect to $C$ evaluated in the centre of the disk is $p$. The task therefore is to evaluate the mapping of $\Omega$ onto $C$ and thence 
the boundary mapping
\begin{equation}
 \label{eq: map}
 \tilde{f} : \Gamma_d \rightarrow \tilde{\Gamma_d}.
\end{equation}
Then, from the above cited lemma, the harmonic measure may be readily evaluated from the boundary mapping and the location of the 
prevertices on the boundary of $C$. The convergence rate of $p$ $-$ with its value characterising 
certain aspects of the Brownian process $-$ helps establish the accuracy of the conformal map.
In the next two subsections we shall consider two methods of 
conformal mapping: first by the Schwarz-Christoffel mapping of the unit disk onto appropriately constructed straight-edge polygons 
(of which the stadium will be a limit); and second by the solution of Symm's integral equation, 
which is applicable in theory to arbitrary geometries.

\subsection{Schwarz-Christoffel mapping} \label{sec: SC-map}
\begin{figure}[hb!]
    \centering
    \begin{subfigure}[h!]{\textwidth}
        \centering
        \includegraphics[scale=0.35]{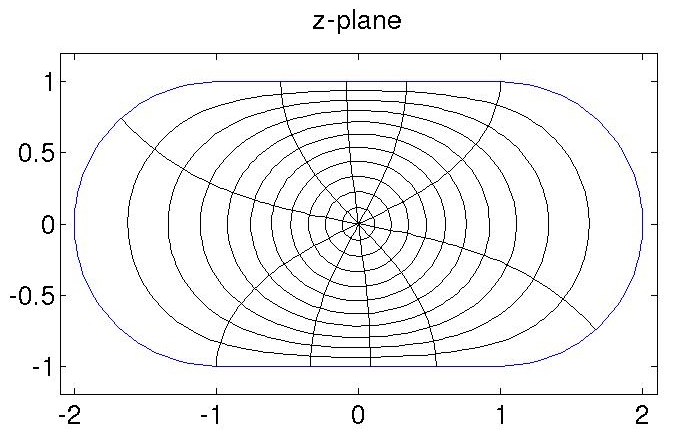}
        \caption{Conformal map, with an estimated error of $10^{-8}$, of the $n=100$ sided polygon approximation to the stadium.}
    \end{subfigure}\\
    ~ 
    \begin{subfigure}[h!]{\textwidth}
        \centering
        \includegraphics[scale=0.35]{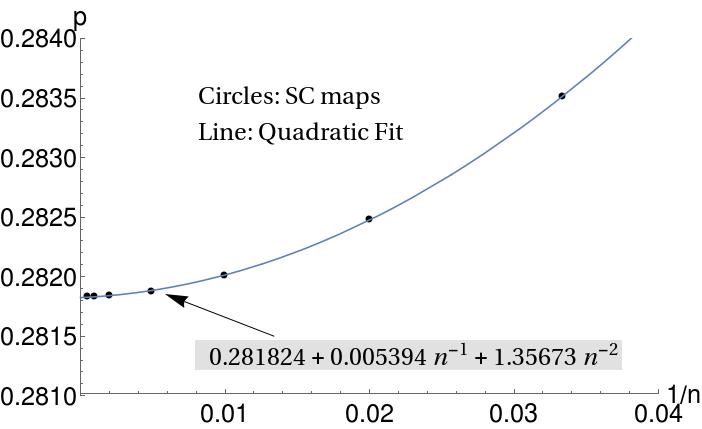}
        \caption{Extrapolation of harmonic measures of the ends evaluated at the centre 
        to the $n \rightarrow \infty$ by a quadratic fit.}
    \end{subfigure}
    \caption{Schwarz-Christoffel maps from the unit disk to $n-$sided polygon approximations to the stadium.}
            \label{fig: SCMap}
\end{figure}
The Schwarz-Christoffel mapping \cite{DT} for straight-edge polygons was, historically, the first formulation of conformal 
transformations of non-trivial geometries. The \textit{Schwarz-Christoffel formula} \cite{DT} describes a conformal 
mapping $f$ from the canonical domain to 
the interior of the given polygon of $n$ sides with interior angles $\alpha_i\pi$
\begin{equation}
 \label{eq: SCmap}
 f(z) = A + C\int^{z}\prod_{k=1}^n\left (1 - \frac{\eta}{z_k}\right) ^{\alpha_k - 1} d\eta,
\end{equation}
where $z_k$ are the prevertices of the polygon on the boundary of the unit disk, the constants $A$ and $C$ fixing 
orientation and scaling of the map. The determination of the prevertices $z_k$ constitute the \textit{parameter problem} 
and is achieved by solving a set of 
equations formed from \eqref{eq: SCmap} that fix the lengths of the each of the sides of the polygon; this set of equations is generally 
solved by some numerical iterative procedure \cite{DT}. Solution of the parameter problem and the evaluation of the maps using 
\eqref{eq: SCmap} is implemented in this work via the Schwarz-Christoffel (SC) Toolbox of MATLAB \cite{Dri}. \\
For $n \le 3$ the conformal mapping, and thence the harmonic measure, may be readily computed due to the freedom in 
choosing the prevertices assured by the Riemann mapping theorem; 
for $n = 4$ rectangular geometries as well the map can be accurately computed using the Schwarz-Christoffel formula \cite{DT, SIAM}. From which the harmonic measure of region's boundary arcs 
may in turn be accurately determined; in particular the harmonic measure of the ends of the rectangle (i.e. its vertical sides) evaluated at the centre is given by the solution of a simple integral equation \cite{SIAM, DT}
\begin{equation}
 \label{eq: SCRect}
 \cot^{-1}{L} = \frac{p\pi}{2} + \arg{K(e^{ip\pi})},
\end{equation}
where $K(k)$ is the complete elliptic integral of the first kind with modulus $k$ and $L$ is the 
ratio of the rectangle's horizontal to vertical side lengths; 
\eqref{eq: SCRect} may be taken to be a relatively exact solution because of the high accuracy with which it can be solved. \\
However for $n-$sided polygons where $n$ is large (e.g. $n > 500$) the crowding of prevertices can lead to 
numerical inaccuracies in the determination of their locations, even using double precision. This may be circumvented by the use of the cross-ratios 
representation implemented in the Toolbox 
\cite[User's Guide]{Dri}: here the crowding induced for large $n$ polygons by the selection of a 
particular choice of the first 3 prevertices (the freedom in doing so being guaranteed by the Riemann mapping theorem) is avoided by 
not computing the prevertices themselves but ratios of 4-tuples of them. In so doing all equivalent (i.e. by a M\"{o}bius transformation) maps 
are represented by a single family, which are given by this set of cross-ratios. The cross-ratios representation may now be converted to the 
standard representation to give the locations of the prevertices; we have found that directly evaluating the disk maps for 
large$-n$ polygons 
leads to non-convergence of the solution. However with the procedure outlined above an estimated error of $10^{-9} - 10^{-7}$ is obtained for 
the conformal maps of all polygons considered. 
From these prevertices the harmonic measure of $\Gamma_d$ may be computed as explained in the paragraph preceeding \eqref{eq: map}.\\
\begin{center}
    \begin{table}[hb!]
    \begin{tabular}{| l | l |}
    \hline
    $n$ & $p = \Delta \theta/\pi$ \\ \hline
     $30$ & 0.28351018 \\ \hline
     $50$ & 0.28247908 \\ \hline
     $100$ &  0.28201121\\ \hline
     $204$ & 0.28187880 \\ \hline
     $500$ & 0.28183922 \\ \hline
     $1000$ &  0.28183244\\ \hline
     $2000$ &  0.2818305\\ \hline
     $n \rightarrow \infty$ & 0.281824 \\
    \hline
    \end{tabular}
     \caption{Harmonic measure $p$ obtained from the Schwarz-Christoffel conformal maps of the unit disk to $n-$sided polygons 
     approximating the stadium. $\Delta \theta$ is defined in \eqref{eq: dtheta}.}
     \label{tab: SCMaps}
    \end{table}
\end{center}
\subsubsection*{Results}
Our strategy then is as follows: the $n-$sided polygon approximation to the stadium is constructed by approximating the 
dome with $(n-2)/4$ straight sides of 
appropriate length that may be inscribed within the unit semi-circle at each end; as $n \rightarrow \infty$ the stadium is realised. 
For instance, a conformal map of the unit disk to the $100-$sided polygon (49 straight edges approximating each semi-circle)
approximation to the stadium is shown in 
Fig. \ref{fig: SCMap}(A), which is visually indistinguishable from the stadium; the ellipses and lines are maps of circles and 
outward drawn radii in $C$. For each value of $n$ the harmonic measure is given by 
\begin{equation}
 \label{eq: dtheta}
 \frac{\Delta \theta}{\pi} = \frac{\theta^{v_1} - \theta^{v_2}}{\pi}, 
\end{equation}
where $\theta^{v_{1,2}}$ are the angles of the prevertices corresponding to $\zeta_0(1)$ and $\zeta_2(-1)$. 
The results of these computations are summarised in Table \ref{tab: SCMaps}.\\
We observe that as $n$ increases we are able to achieve at most 4-5 digits of consistency amongst the results for the last four $n$ values 
i.e. $p \approx 0.2818(3)$; thus $10^{-5}$ may be taken to as the absolute error in the conformal transformation of the stadium's 
region $\Omega$ as approximated by the largest polygon.
For the harmonic measure, however, we may straightforwardly extrapolate the data, by a quadratic fit, to estimate 
the $n \rightarrow \infty$ limit; this extrapolation is shown in Fig. \ref{fig: SCMap}(B), which gives $p \approx 0.281824$. We will see 
in the next section that the fifth digit from the extrapolation is also correct.
\subsection{Symm's integral equation} \label{sec: Symms}
The conformal mapping problem of a region is equivalent to determining a harmonic function $g(z)$ (and its harmonic conjugate obeying 
the Cauchy-Riemann equations) that satisfy logarithmic boundary conditions on that region \cite{Symm}. 
$g(z)$ when represented as a logarithmic potential (which preserves harmonicity)
\begin{equation}
 \label{eq: Symms1}
 g(z) = \int_{\Gamma} \log{|z - \eta|}\sigma (\eta) |d\eta|,
\end{equation}
reduces the problem further to that of determining the \textit{source density} $\sigma (\eta)$ such that a set of equations, now called Symm's equations, 
are satisfied \cite{Symm}. In our representation of the boundary $\Gamma$,  parametrised by $t$ for each of four arcs in 
\eqref{eq: Bunimovich}, Symm's equations for the mapping of the interior $\Omega$ onto the unit disk take the form 
\cite{Symm, Hough}
\begin{eqnarray}
\label{eq: Symms2}
\sum_{k=0}^{3}\int_{-1}^{1}\sigma_k(t)\log{|z - \zeta_k(t)|}dt &=& \log{|z|}, \hspace{4em} z \in \Gamma, \nonumber \\
\sum_{k=0}^{3}\int_{-1}^{1}\sigma_k(t) &=& 1.
\end{eqnarray}
A source density $\sigma_k(t)$ has now been associated with each arc $\zeta_k$. 
The conformal map to the unit disk for a point $z \in \Omega$ (which may also be extended to the boundary 
$\Gamma$ appropriately) is then given by \cite{Symm, Hough}
\begin{equation}
\label{eq: P1}
 f(z) = z\exp{(-P(z))},
\end{equation}
where $P(z)$ is obtained from the solutions to \eqref{eq: Symms2} as
\begin{equation}
 \label{eq: P2}
 P(z) \defeq \sum_{n=0}^{3}\int_{-1}^{1}\sigma_k(t)\log{(z - \zeta_k(t))}dt.
\end{equation}

The method we adopt for solving for the source densities $\sigma_k(t)$ is that of 
collocation by Chebyshev polynomials as introduced in Ref. \cite{Hough}; the singularities in each 
$\sigma_k$ at its end points
are approximated, generically and independent of the precise nature of the singularity, by the usual Chebyshev weights 
$w = \frac{1}{\sqrt{1-t^2}}$. This makes it a 
relatively simple and practical technique for conformal mapping; however, as Hough et al. point out 
\cite[Theorem 2.1]{Hough}, for smooth boundaries 
(which is by and large also our case due to the absence of any sharp corners) the convergence rate is rather slow. 
As expected we shall 
see that the convergence rate of the maps for the rectangular geometries is much more superior. Nevertheless the numerical accuracy of the 
solution for the harmonic measure of $\Gamma_d$ fares better than the two techniques implemented in the previous sections, thereby assuring 
a better conformal map.\\
We summarise this technique adapted from Ref. \cite{Hough}: Within this collocation scheme and approximation of the 
singularities the source density on the $k^{\textrm{th}}$ arc is expressed as
\begin{equation}
 \label{eq: SymmsDensity}
\sigma_k(t) \defeq \frac{\phi_k(t)}{\sqrt{1 - t^2}},
\end{equation}
where the analytic functions $\phi_k(t)$ are approximated by a collocation expansion with coefficients $\phi_{kn}$
\begin{equation}
 \label{eq: SymmsCollocation}
 \tilde{\phi_k}(t) \defeq \sum_{n=0}^{\nu_k}\phi_{kn}T_n(t),
\end{equation}
where $\nu_k$ denotes the number of collocation points on the $k^{\textrm{th}}$ arc and $T_n(t)$ is the 
$n^{\textrm{th}}$ order Chebyshev 
polynomial of the first kind. Using \eqref{eq: SymmsDensity} and \eqref{eq: SymmsCollocation} in \eqref{eq: Symms2} 
results in a system of linear equations in $\{\phi_{kn}\}$
\begin{eqnarray}
 \label{eq: SymmsEquations}
 \sum_{k=0}^3\sum_{n=0}^{\nu_k}C_{jmkn}\phi_{kn} &=& \mu_{jm}; \hspace*{3em} j=0,1,2,3 \hspace*{2em} m=0,1,2,\ldots d_j , \nonumber \\
 \pi\sum_{k=0}^3\phi_{k0} &=& 1.
\end{eqnarray}
The coefficients $C_{jmkn}$ and constants $\mu_{jm}$ are defined by 
\begin{eqnarray}
 \label{eq: SymmsDefinitions}
 C_{jmkn} &\defeq& \int_{-1}^{1}\frac{T_n(t)}{\sqrt{1-t^2}}\log{|\zeta_j(\tau_m^{(j)}) - \zeta_k(t)|}dt, \nonumber \\
 \mu_{jm} &\defeq& \log{|\zeta_j(\tau_m^{(j)})|},
\end{eqnarray}
where $\tau_m^{(j)} \defeq \cos{\left[\cfrac{(2m + 1)\pi}{2\nu_j + 2}\right]}$ parametrises the $m^{\textrm{th}}$ collocation point on the 
$j^{\textrm{th}}$ arc; we take $d_j = \nu_j = \nu \gg 1$. This gives an overdetermined system of $4\nu + 5$ 
equations in $4\nu + 4$ coefficients. The numerically computed residuals for such a system of equations for the largest 
$\nu$ are found to be between $10^{-15} - 10^{-10}$.\\
Once the $\{\phi_{kn}\}$ are computed the conformal map $f(z)$ to the unit disk may be approximated with 
\eqref{eq: P1} and \eqref{eq: P2}; the corresponding angles $\theta_k(t)$ on the unit disk 
of the $k^{\textrm{th}}$ boundary arc may in turn be approximated as \cite{Symm,Hough} 
\begin{eqnarray}
 \label{eq: SymmsFinal}
\theta_k(t) = \theta_k(-1) + 2\pi\phi_{k0}(\pi - \arccos{(t)}) - 2\pi\sqrt{1-t^2}\sum_{n=1}^{\nu}\frac{\phi_{kn}U_{n-1}(t)}{n},
\end{eqnarray}
where $U_n(t)$ is the Chebyshev polynomial of the second kind. This gives the location of the prevertices, and the conformal map of the 
region $\Omega$ and the boundary $\Gamma$ is complete. Finally, using \eqref{eq: SymmsFinal}, the harmonic measure of 
the segment $\Gamma_d$ is given by
\begin{equation}
 \label{eq: HM}
 p \equiv \frac{\theta_1(1) - \theta_1(-1)}{2\pi} + \frac{\theta_3(1) - \theta_3(-1)}{2\pi} = \pi\phi_{10} + \pi\phi_{30}.
\end{equation}
We note that by the symmetry of the boundary $\Gamma$, $\phi_{10} = \phi_{30}$.

\subsubsection*{Coefficients evaluation}

The main computational resources are spent in evaluating the integrals in \eqref{eq: SymmsDefinitions}. 
Some simplification is possible because of the 
symmetry of the domain $\Omega$. In particular, the following relations may be seen to hold among the coefficients and 
constants in \eqref{eq: SymmsDefinitions}
\begin{eqnarray}
 \label{eq: Symmetry}
 C_{jmkn} &=& C_{(j+2)(\textrm{mod} 4),m,k,(n+2)(\textrm{mod} 4)}, \nonumber \\
 C_{0m3n} &=& (-1)^{n}C_{0,\nu-m, 1, n}, \nonumber \\
 C_{3m0n} &=& (-1)^{n}C_{1,\nu-m, 0, n}, \nonumber \\
 \mu_{jm} &=& \mu_{(j+2)(\textrm{mod}4), m}.
\end{eqnarray}
These symmetries reduce the number of integrals to be computed by over a factor of 2.\\
A main of source of difficulty in evaluating the integrals in \eqref{eq: SymmsDefinitions} arises when $j=k$ and a 
logarithmic singularity is encountered in the numerics. This may be circumvented by noting that either an exact expression
is available for the integral \cite{Prudnikov, Hough} or a standard singularity removal procedure \cite{Hough} may be 
employed for better convergence. \\
Consider first $j=k=0$: then $C_{0m0n}$ in \eqref{eq: SymmsDefinitions} is given by
\begin{equation}
\label{eq: C0m0n}
C_{0m0n} = \int_{-1}^{1}\frac{T_n(t)}{\sqrt{1-t^2}}\log{|\tau_m^{(0)} - t|} dt;
\end{equation}
this may be exactly evaluated \cite{Prudnikov, Hough} to give
\begin{equation}
  C_{0m0n}=\begin{cases}
    -\frac{\pi}{n}\cos{(n\arccos{(\tau_m^{(0)})})}, & \text{$n>0$}.\\
    \pi\log{\left(\frac{L}{2}\right)}, & \text{$n=0$}.
  \end{cases}
\end{equation}
Next consider $j=k=1$: then $C_{1m1n}$ in \eqref{eq: SymmsDefinitions} is given by 
\begin{equation}
\label{eq: C1m1n}
C_{1m1n} = \int_{-1}^{1}\frac{T_n(t)}{\sqrt{1-t^2}}\log{|\zeta_1(\tau_m^{(1)}) - \zeta_1(t)|} dt;                              
\end{equation}
substituting $t = \cos{\theta}$ and defining $\alpha_m \defeq \frac{(2m+1)\pi}{2\nu+2}$, the above integral can be 
reexpressed as 
\begin{equation}
 \label{eq: C11}
 C_{1m1n} = \int_0^{\pi}\cos{(n\theta)}\log{\left[ \frac{2\sin{((\alpha_m-\theta)/2)}}{\alpha_m-\theta}\right]} 
 + \int_0^{\pi}\cos{(n\theta)}\log{|\alpha_m - \theta|}d\theta.
\end{equation}
The first integral on the right hand side may be evaluated numerically; for $n \neq 0$ the second integral may be 
simplified, after some algebra and integration by parts, to give 
\begin{eqnarray}
 \label{eq: C11n}
 \int_0^{\pi}\cos{(n\theta)}\log{|\alpha_m - \theta|}d\theta &=& 
 \frac{\sin{(2n\alpha_m)}}{2n}\Big[ \textrm{ci}(n\alpha_m) - \textrm{si}(n\alpha_m) \nonumber \\
 &-& \textrm{ci}(n\pi - n\alpha_m) - \textrm{si}(n\pi - n\alpha_m)\Big],
\end{eqnarray}
where $\textrm{ci}(z) = -\int_z^{\infty}\frac{\cos{t}}{t}dt$ is the cosine integral and 
$\textrm{si}(z) = \int_0^{z}\frac{\sin{t}}{t}dt$ is the sine integral.\\
For the special case $n=0$, \eqref{eq: C11} may be expressed in terms of the Clausen function 
$\textrm{Cl}_2(z) = -\int_0^z\log{(2\sin(\theta/2))}$
\begin{equation}
 \label{eq: C110}
 C_{1m10} = -\textrm{Cl}_2(\alpha_m) - \textrm{Cl}_2(\pi - \alpha_m).
\end{equation}
All the terms in \eqref{eq: C11n} and \eqref{eq: C110} may be efficiently evaluated to give good accuracy for $C_{1m1n}$. 
The other two singular coefficients $C_{2m2n}$ and $C_{3m3n}$ are, from \eqref{eq: Symmetry}, equal to 
$C_{0m0n}$ and $C_{1m1n}$ respectively. All the coefficients $C_{jmkn}$ may now be accurately evaluated for 
forming the system of linear equations \ref{eq: SymmsEquations}. We employ the GSL library \cite{GSL} for integration and 
solving the system of linear equations, and insist on 
the estimated absolute error for the integrals to be \textit{at least} $10^{-8}$.

\subsubsection*{Results}
\begin{figure}[t!]
    \centering
    \begin{subfigure}[h!]{\textwidth}
        \centering
        \includegraphics[scale=0.35]{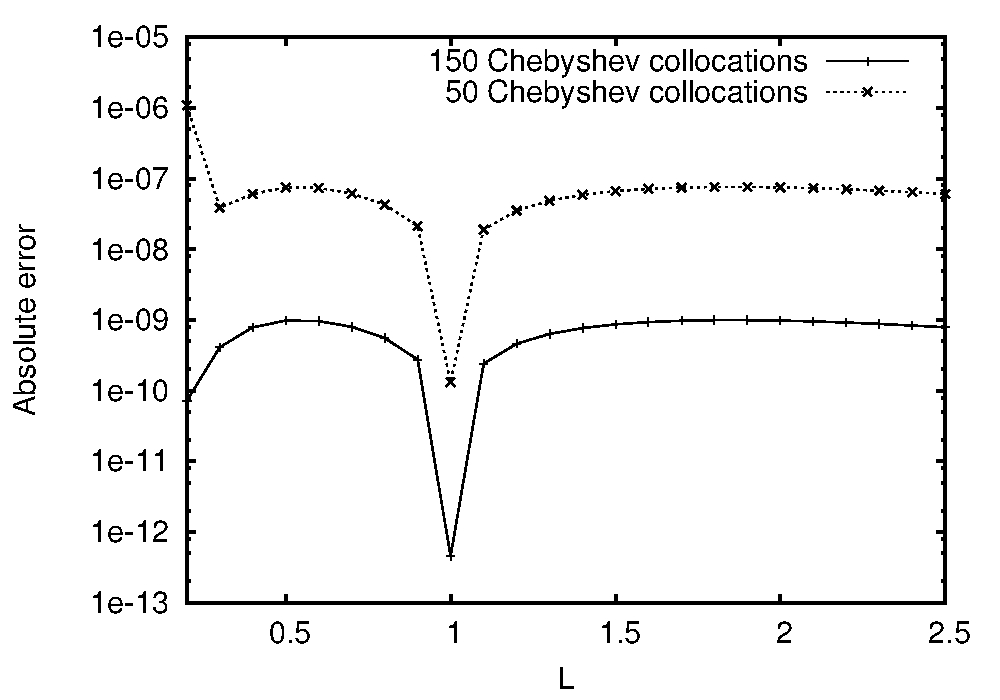}
        \caption{Difference of ends' harmonic measures for the $2L \times 2$ rectangles between the solutions from the 
        ``exact'' Schwarz-Christoffel formula \eqref{eq: SCRect} and Symm's integral equations \eqref{eq: SymmsEquations} for 
        two different collocations in each arc.}
    \end{subfigure}\\
    ~ 
    \begin{subfigure}[h!]{\textwidth}
        \centering
        \includegraphics[scale=0.35]{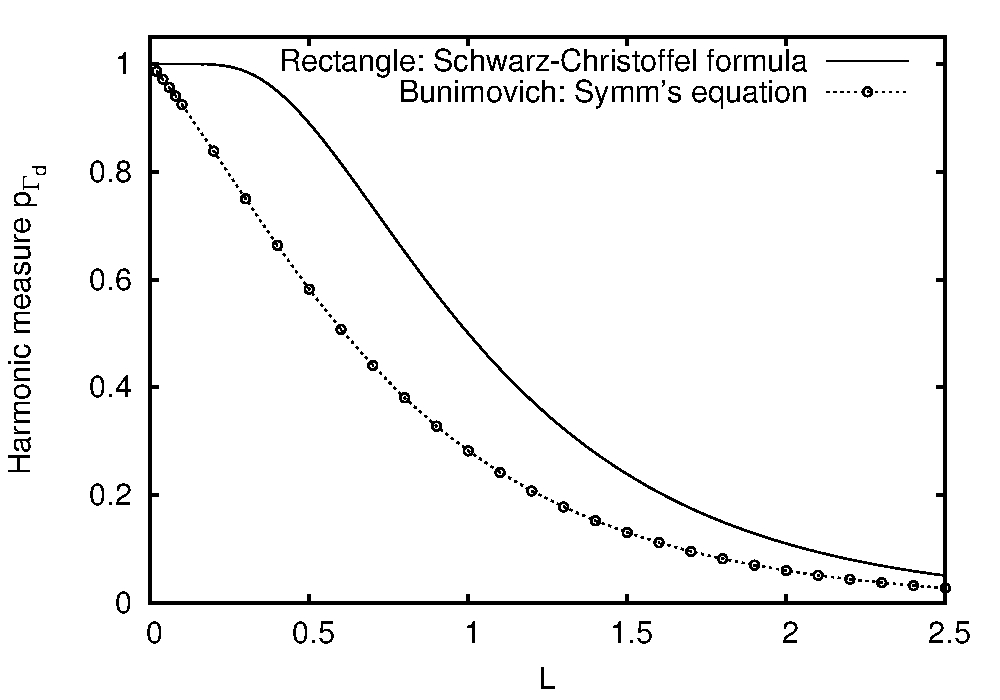}
        \caption{Harmonic measure of the ends in the Bunimovich stadium (200 collocations in Symm's equations 
        \eqref{eq: SymmsEquations}) versus the rectangular geometry (Schwarz-Christoffel formula \eqref{eq: SCRect}) for 
        various lengths $L$.}
    \end{subfigure}
    \caption{Solutions of Symm's integral equation for the Bunimovich stadium and the rectangular geometry.}
            \label{fig: StadiumRectangular}
\end{figure}

As a testbed for the above procedure we first construct conformal maps for rectangular geometries, and compute $p$ in order to compare 
with its highly accurate solution \eqref{eq: SCRect}; 
the symmetries of the rectangle still retain the validity of \eqref{eq: Symmetry}, and the handling of the 
numerical singularities remains essentially unchanged. We plot in Fig. \ref{fig: StadiumRectangular}(A) the absolute 
error between the ``exact'' results for the harmonic measure of the rectangle's ends \eqref{eq: SCRect} and 
the values obtained from solving Symm's equations \eqref{eq: SymmsEquations}. 
Recall that from the latter solutions $\phi_{kn}$ we use \eqref{eq: HM} to compute $p$. We note from Fig. \ref{fig: StadiumRectangular}(A) 
that even for a modest number of collocations the absolute error in $p$ decreases very rapidly for a range of side lengths $L$; Symm's 
equations indeed give highly precise conformal maps for the rectangular geometries with the expenditure of little effort.
We will see that the situation for the Bunimovich stadium is very different in that the convergence is 
limited even for large $\nu$ values, quite likely due to the smoothness of the 
boundary $\Gamma$ for the stadium and the artificial introduction of singularities in the source densities 
$\sigma(t)$ in \eqref{eq: SymmsDensity}.\\
Consider now the Bunimovich stadium for various lengths $L$ of the straight edges as shown in Fig. \ref{fig: stadium}. The conformal 
transformations are computed using Symm's equations \eqref{eq: SymmsEquations}
with $\nu = 200$ collocation points in each arc. The harmonic 
measure $p$ of the domes $\Gamma_d$ is extracted as was done for the rectangular geometry. 
It is clear that for a given $L$, 
we expect $p^{\textrm{rectangle}} > p^{\textrm{stadium}}$ for the ends; the results for the harmonic measure are shown in 
Fig. \ref{fig: StadiumRectangular}(B) for the two geometries, and this expectation is borne out. \\
However for $L \ll 1$ 
it is interesting to contrast the behaviour of the harmonic measure in the two cases i.e. the presence of a plateau in 
the rectangle versus the sudden drop in stadium; the Brownian particle diffuses towards the horizontal edges in the 
stadium much more rapidly than in the rectangle. It would be interesting, as future study, to analytically 
corroborate this sudden drop in the harmonic measure 
as a limiting case of the circle being infinitesimally deformed by horizontal straight edges, which corresponds to 
$L \ll 1$ for the stadium.\\
Now we focus on the length $L=1$, which was studied in Sec. \ref{sec: MC} with Monte Carlo and in Sec. \ref{sec: SC-map} with 
Schwarz-Christoffel transformations, in order to push the accuracy of map and simultaneously the measure $p$. Symm's equations were solved for a 
range of collocation points $\nu$ from $64-1200$; from the computed conformal maps for each case the estimate of $p$ is 
listed in Table \ref{tab: Symms}. We see that, unlike in the rectangular case, the convergence rate (here measured by 
agreement between results from succeeding $\nu$ values) is not too rapid. \\
\begin{center}
    \begin{table}[t!]
    \begin{tabular}{| l | l |}
    \hline
    $\nu$ & $p = \pi\phi_{10} + \pi\phi_{30}$ \\ \hline
     $64$ &   ${\bf 0.281}76556$\\ \hline
     $100$ &  ${\bf 0.2818}0170$\\ \hline
     $128$ &  ${\bf 0.2818}1209$\\ \hline
     $256$ &  ${\bf 0.28182}502$\\ \hline
     $300$ &  ${\bf 0.28182}628$\\ \hline
     $350$ &  ${\bf 0.28182}718$\\ \hline
     $500$ &  ${\bf 0.28182}850$\\ \hline
     $512$ &  ${\bf 0.28182}856$\\ \hline
     $800$ &  ${\bf 0.281829}30$\\ \hline
     $1000$ & ${\bf 0.281829}50$\\ \hline
     $1200$ & ${\bf 0.281829}60$\\ \hline
    \end{tabular}
     \caption{Harmonic measures obtained from the solutions of Symm's equations of the Bunimovich stadium 
     for various number of Chebyshev collocation points $\nu$.}
     \label{tab: Symms}
    \end{table}
\end{center}
At values as small as $\nu=64$ (which takes the 
program a few seconds to complete), the value of $p$ agrees with the largest Monte Carlo runs (which took a few hours to complete). As 
the number of collocations are increased the digits that agree with the succeeding result are shown in bold in the table. 
From collocations up to and including $\nu = 1200$, we are thus able to assert six digits for the harmonic measure
\begin{equation}
 \label{eq: Final}
 p = 0.281829\ldots
\end{equation}
This exemplifies the point earlier made, as well in Ref. \cite{Hough}, that the convergence rate of the map obtained with 
this method is of order $\nu^{-2}$ for smooth boundaries; therefore to get an accuracy of 10 digits we expect that 
about $10^5$ collocation points would be required.
Note also that the fifth digit from the naive quadratic extrapolation of the results from the Schwarz-Christoffel maps 
in Fig. \ref{fig: SCMap}(B) agrees with \eqref{eq: Final}. We have thus obtained a numerical conformal map of the stadium, using this 
readily programmable algorithm of Hough et al. \cite{Hough}, to within an error of $10^{-7}$.\\
We are confident that the conformal map and the accuracy obtained in 
\eqref{eq: Final} may be bettered by the use of other collocation schemes that treat the singularities 
in the source densities individually or by a different numerical conformal mapping technique altogether. Of course, in order to merely 
obtain higher precision for $p$, one might consider using direct finite-element solvers of the Laplacian \eqref{eq: Poisson}.
\section{Conclusions} \label{sec: Conclusions}
We have studied the problem of a Brownian particle in a Bunimovich stadium $\Omega$ with Dirichlet boundary conditions.
Firstly a stochastic simulation of the process was employed. Further a conformal mapping of the region onto the unit disk was also undertaken; 
this latter approach utilises a deterministic reformulation of the problem, which enables for a more precise analysis.
By conformally mapping the stadium to the unit disk, the harmonic measures of the boundary arcs are readily evaluated. \\
For the conformal transformations, two methods were undertaken: 
Schwarz-Christoffel mapping of $n-$sided polygons in the large$-n$ limit, and the solution of Symm's integral 
equations via a Chebyshev weighting of the solutions; using which we are able to achieve $5$ and $6$ digits of accuracy for the 
harmonic measure, which 
can certainly be improved by other collocation schemes or different numerical conformal mapping algorithms. The convergence rate of 
the harmonic measure gives a gauge on the accuracy of the stadium's conformal map: in particular, using Chebyshev weighted solutions Symm's equations we obtained a conformal map of the stadium to within an absolute error of $10^{-7}$. Interesting 
differences in the behaviour of the harmonic measure for the limit $L \ll 1$ are exhibited between the rectangular and 
stadium geometries, which we believe might be corroborated analytically for infinitesimally deformed circles.

\begin{acknowledgement}
We thank H. Monien and T.A.~Driscoll for helpful discussions and correspondence.\end{acknowledgement}

\end{document}